\documentclass{article}
\usepackage{amsthm,amsmath,amssymb}
\RequirePackage[numbers]{natbib}
\RequirePackage[colorlinks,citecolor=blue,urlcolor=blue]{hyperref}

\def\bE{\mathbb E}
\def\bE{\mathbb E}
\def\R{\mathbb R}

\newtheorem{Theorem}{Theorem}

\newtheorem{lemma}[Theorem]{Lemma}
\newtheorem{Corollary}[Theorem]{Corollary}
\newtheorem{proposition}[Theorem]{Proposition}
\newtheorem{Assumption}{Assumption}
\newtheorem{remark}{Remark}

\def\rank{{\rm rank}}

\def\keywords{\vspace{.5em}
{\textit{Keywords}:\,\relax%
}}
\def\endkeywords{\par}

\def\subjclass{\vspace{.5em}
{\textit{AMS 2000 subject classification}:\,\relax%
}}
\def\endsubjclass{\par}

\begin{document}

\title{Matrix completion by singular value thresholding: sharp bounds}
\author{Olga Klopp\\  \\ CREST and MODAL'X, University Paris Ouest}
\maketitle


\begin{abstract}
We consider the matrix completion problem where the aim is to estimate a large data matrix for which only a  relatively small random subset of its entries is observed. Quite popular approaches to matrix completion problem are iterative thresholding methods. In spite of their empirical success, the theoretical guarantees of such iterative thresholding methods are poorly understood. The goal of this paper is to provide strong theoretical guarantees, similar to those obtained for nuclear-norm penalization methods and one step thresholding methods,  for an iterative thresholding algorithm  which is a modification of the softImpute algorithm. An important consequence of our result is the exact minimax optimal rates of convergence for matrix completion problem which were known until know only up to a logarithmic factor. 
\end{abstract}
 
%
\endkeywords

\endsubjclass


\section{Introduction}\label{introduction}
 Suppose that we observe a small subset of entries of a large data matrix. The problem of inferring the many missing entries from this small set of observations is known as the matrix completion problem. This problem has attracted considerable attention in the past five years. 
  The first works \cite{Candes_Recht,Candes_Tao, Candes_Plan, Gross, Recht} introduce  nuclear-norm minimization method. A different approach, called OPTISPACE has been proposed in \cite{Keshavan_Montanari_Oh,Keshavan}. More recently, a method based on  max-norm minimization was studied in \cite{Cai_Zhou,Foygel_Serebro}. Other methods include, for example,  GROUSE (Grassmannian Rank-One Update Subspace Estimation) \cite{Balzano} and orthogonal rank-one matrix pursuit \cite{Wang}. 
  
 A quite popular direction in the matrix completion literature are the thresholding methods  which can be divided in two groups: one-step thresholding methods and iterative thresholding methods. Strong theoretical guarantees were obtained for one-step thresholding procedures. For example, Koltchinskii et al in \cite{Koltchinskii_Lounici_Tsybakov} introduce a soft-thresholding method and show that it is minimax optimal up to a logarithmic factor. In \cite{ Klopp_rank} Klopp consider a hard thresholding proceedure. Chatterjee \cite{ Chatterjee_mc} propose an universal singular value thresholding that can be applied to a large number of matrix estimation problems, including matrix completion. Despite strong theoretical guarantees these one-step thresholding methods has two important drawbacks: they show poor behavior in practice and only work under the uniform sampling  distribution which is not realistic in many practical situations.

  Much better practical performances have been shown by iterative thresholding methods. For example, in \cite{Cai_Candes_Shen}, Cai et al  propose a first-order 
  singular value thresholding
  algorithm SVT  which  approximately
  solves the nuclear norm minimization problem.  In \cite{Mazumder_Hastie_Tibshirani}, Mazmuder et al introduce \texttt{softImpute} algorithm. \texttt{softImpute} produces a sequence of solutions that converges to a solution of   the nuclear norm regularized least-squares problem when the number of iterations goes to infinity. These iterative thresholding  algorithms are simple to implement,
 scale to relatively large matrices and in practice achieve  competitive errors compared to
 the  state-of-the-art algorithms. More recently Dhanjal et al \cite{Dhanjal} propose an improvement for  the \texttt{softImpute} algorithm using randomized SVDs along with a novel updating method. This improvement allows to bypass the bottleneck in the algorithm which consists in the use of the
 singular value decomposition  of a large matrix at each iteration. 
 
 The majority of existing algorithms for matrix completion are batch methods, that is, they operate on the full data matrix. However in some applications such as recommendation systems or localization in sensor networks  we observe a sequence of data matrix $M_1,\dots ,M_T$ reviled sequentially where from $M_t$ to $M_{t+1}$ we add new observations.  In such situations
  the predictive rule should be refined incrementally. One advantage of  iterative thresholding algorithms is that they can be adapted to such sequential learning, see for example \cite{Dhanjal}.
  
  In spite of their empirical success, the theoretical guarantees of such iterative thresholding methods are poorly understood. The goal of this paper is to provide strong theoretical guarantees, similar to those obtained for nuclear-norm penalization methods (see, for example \cite{Negahban_Wainwright,Klopp_general}) and one step thresholding methods (see \cite{Koltchinskii_Lounici_Tsybakov,Klopp_rank,Chatterjee_mc}) for a modification of the \texttt{softImpute} algorithm.
  \subsection{Contributions and Related Work}
   The contributions of the present paper to the theoretical study of the modified \texttt{softImpute} algorithm are multifaceted. In Section \ref{section_upper_bounds} we prove an upper bound on the estimation error of the output $\hat M$ of our algorithm. Let $M_0\in \mathbb{R}^{m_1\times m_2}$ be the unknown matrix of interest. Suppose, for simplicity, that each entry is observed with the same probability $p$, then  we prove the following upper bound on the estimation error of $\hat M$
   \begin{equation} \label{upper_bound}
   \begin{split}
   \dfrac{\Vert \hat M-M_0\Vert_{2}^{2}}{m_1m_2}&\lesssim \dfrac{\rank(M_0)}{p\min(m_1,m_2)}.
   \end{split}
   \end{equation}
Here the symbol $\lesssim$ means that the inequality holds up to a multiplicative numerical constant. To the best of our knowledge, the upper bound on the estimation error given by \eqref{upper_bound} is strictly better than all  upper bounds  available in matrix completion literature.

 For instance, for the same setting, Chatterjee in  \cite{Chatterjee_mc} obtains the following larger bound
 \begin{equation*} 
    \begin{split}
    \dfrac{\Vert \hat M-M_0\Vert_{2}^{2}}{m_1m_2}&\lesssim \sqrt{\dfrac{\rank(M_0)}{p\min(m_1,m_2)}}.
    \end{split}
    \end{equation*}
  On the other hand, \cite{Koltchinskii_Lounici_Tsybakov,Negahban_Wainwright,Klopp_general}, among some other papers, consider a slightly different  setting where the matrix completion problem is viewed as a particular case of the trace regression model. In this setting the number of observations $n$ is fixed.  The drawback here is that in this model each entry  can be observed multiple times which is not the case in a large number of practical situations.  We consider a  different setting where each entry can be observed at most once (see Section \ref{model_sampling}). However, it is easy to see that these two settings are closely related if we put $n=pm_1m_2$. Comparing to \eqref{upper_bound}, the bounds obtained in \cite{Koltchinskii_Lounici_Tsybakov,Negahban_Wainwright,Klopp_general} have an additional $\log(d_1+d_2)$ factor.
  
  Koltchinskii et al in \cite{Koltchinskii_Lounici_Tsybakov} obtained lower bounds for the estimation error without this additional $\log(d_1+d_2)$ factor. So our result answer the important theoretical question what is the exact minimax rate of convergence for matrix completion problem. As the lower bound in \cite{Koltchinskii_Lounici_Tsybakov} is obtained for a different setting, in Section \ref{lower_bounds} we adapt their proof to our setting, showing that the minimax rate of convergence for matrix completion problem is given by \eqref{upper_bound} and that the estimator produced by our algorithm is minimax optimal. Note that our techniques can be adapted to the setting considered in \cite{Koltchinskii_Lounici_Tsybakov,Negahban_Wainwright,Klopp_general} and lead to an upper bound without the additional $\log(d_1+d_2)$ factor in this setting also.
  
  Another important point is that a large part of matrix completion literature consider uniform sampling at random setting where each entry is observed with the same probability $p$. In many applications, such as recommendation systems, this assumption is not realistic. The theoretical analysis in the present paper is carried out for quite general sampling distributions and show that our iterative thresholding algorithm has good performances in such situations.
   Finally our results give theoretical insights for the chose of the parameters in the modified \texttt{softImpute} algorithm.
   
   \subsection{Organisation of the paper}
   The remainder of this paper is organized as follows. In Section \ref{model_sampling} we introduce our model and the assumptions on the sampling scheme. For the reader's convenience, we  collect notation which we use throughout the paper in Section \ref{notation}. In Section \ref{section_algorithm} we present a modification of the \texttt{softImpute} algorithm for matrix completion. The upper bounds on the estimation error are derived in Section \ref{section_upper_bounds}. Finally the lower bounds are obtained in  Section \ref{lower_bounds} and the Appendix contains the proofs.

\section{Preliminaries}\label{Preliminaries}
\subsection{Model and sampling scheme}\label{model_sampling}

Suppose that we observe a relatively small number of entries of a data matrix \begin{equation}\label{model}
X=M_0+E.
\end{equation}
 Here $M_0=(m_{ij})\in \mathbb{R}^{m_{1}\times m_{2}}$ is the unknown matrix of interest and $E=(\xi_{ij})\in \mathbb{R}^{m_{1}\times m_{2}}$ is the matrix containing the noise. We assume that the noise variables $\xi_{ij}$ are independent, zero mean and bounded:
 \begin{Assumption}\label{noise_bounded}
$\bE(\xi_{ij})=0$, $\bE(\xi_{ij}^{2})=\sigma^{2}$ and there exists a positive constant  $b>0$ such that
\begin{equation*} 
 \underset{i,j}{\max}\left \vert\xi_{ij}\right \vert\leq b.
 \end{equation*}
 \end{Assumption} 
We suppose that each entry of $X$ is observed independently of the other entries. For the entry $(i,j)\in [m_1]\times[m_2]$, we  denote the probability to be observed by $\pi_{ij}$. Let $\eta_{ij}$ be the independent Bernoulli variables with parameters $\pi_{ij}$  and $y_{ij}=\eta_{ij}\left (m_{ij}+\xi_{ij}\right )$. Then, $Y=(y_{ij})$ is the matrix containing our observations.
We denote by $\Omega$ the random set of observed indices.

In the simplest situation each coefficient is observed with the same probability, i.e. for every $(i,j)\in [m_1]\times[m_2]$, $\pi_{ij}=p$.
  Unfortunately, such an assumption on the sampling distribution is not realistic in many practical applications. In the present paper, we consider general sampling model. 
We suppose that each coefficient is observed with a positive probability:
 \begin{Assumption}\label{sampling_1} There exists  $p>0$ such that for any $(i,j)\in \{1,\dots,m_1\}\times \{1,\dots,m_2\}$
$$\pi_{ij}\geq p.$$
\end{Assumption}

For any $A=(A_{ij})\in \mathbf{R}^{m_1\times m_2}$ we define the weighted by $\pi_{ij}$ Frobenius norm of $A$
\begin{equation*}
\Vert A\Vert _{L_2(\Pi)}^{2}=\sum_{(i,j)}\pi_{ij}A^2_{ij}.
\end{equation*}
Assumption \ref{sampling_1} implies that \begin{equation}\label{ass1}
 \Vert A\Vert^{2} _{L_2(\Pi)}\geq p^{-1}\Vert A\Vert^{2} _{2}.
 \end{equation}  

We denote the column and row marginals by $$\pi_{\cdot j}=\underset{i=1}{\overset{m_1}{\Sigma}}\pi_{ij}\qquad\text{and}\qquad  \pi_{i\cdot }=\underset{j=1}{\overset{m_2}{\Sigma}}\pi_{ij}.$$
 Suppose that we know an upper bound $L$ on it's maximum:
\begin{equation}\label{upper_bound_marginals}
\underset{i,j}{\max}\left (\pi_{\cdot j},\pi_{i \cdot}\right )\leq L.
\end{equation}
Note that we can easily get an estimation on this upper bound using the empirical frequencies
\begin{equation*}
\hat\pi_{\cdot j}=\dfrac{\sum_{i=1}^{m_1}\eta_{ij}}{\sum_{(i,j)}\eta_{ij}}\qquad\text{and}\qquad\hat\pi_{i\cdot }=\dfrac{\sum_{j=1}^{m_2}\eta_{ij}}{\sum_{(i,j)}\eta_{ij}}.
\end{equation*}


\subsection{Notation}\label{notation}
We provide a brief summary of the notation used throughout this paper. Let $A,B$ be matrices in $\mathbb{R}^{m_{1}\times m_{2}}$.
\begin{itemize}
\item For a matrix $A$, $A_{ij}$ is its $ (i, j)$−th entry.
\item We denote by $S_{\lambda}(W)\equiv UD_{\lambda}V'$ the \textit{soft-thresholding} operator where $D_{\lambda}=\text{diag}\left [(d_1-\lambda)_{+},\dots,(d_r-\lambda)_{+}\right ]$, $UDV'$ is the SVD of $W$, $D=\text{diag}\left [d_1,\dots,d_r\right ]$ and $t_{+}=\max(t,0)$.
\item For any set $I$, $\vert I\vert$ denotes its cardinal and $\bar I$ its complement. Let $a\vee b=\max(a,b)$ and $a\wedge b=\min(a,b)$.
\item For two matrices $A,B\in \mathbb{R}^{m_1\times m_2}$ we define the \textit{scalar product}
$$\langle A,B\rangle =\mathrm{tr}(A^{T}B).$$
\item We denote by $\Vert A\Vert_{2}$ the usual $l_2-$norm.  Additionally, we use the following matrix norms: $\Vert A\Vert_{*}$ is the nuclear norm (the sum of singular values), $\Vert A\Vert$ is the operator norm (the largest singular value), $\Vert A\Vert_{\infty}$ is the largest absolute value of the entries:  
$$\left\Vert A\right\Vert_{\infty}=\underset{i,j}{\max}\mid A_{ij}\mid.$$ 
\item  $\pi_{ij}$ is the probability to observe the $(i,j)$-th element. For $j=1\dots m_2$, $\pi_{\cdot j}=\underset{i=1}{\overset{m_1}{\Sigma}}\pi_{ij}$  and for $i=1\dots m_1$,  $\pi_{i\cdot }=\underset{j=1}{\overset{m_2}{\Sigma}}\pi_{ij}$. We have that 
\begin{equation*}
\underset{i,j}{\max}\left (\pi_{\cdot j},\pi_{i \cdot}\right )\leq L.
\end{equation*}
\item Let $M=\max(m_1,m_2)$, $m=\min(m_1,m_2)$ and $d=m_1+m_2$.
\item  Let $I\subset \{1,\dots m_1\}\times\{1,\dots m_2\}$ be a subset of indices. Given a matrix $A=(A_{ij})$, we define its restriction on $I$, $A_{I}$, in the following way: $\left (A_{I}\right )_{ij}=A_{ij}$ if $(ij)\in I$ and $\left (A_{I}\right )_{ij}=0$ if not. 

\item We denote $\Vert A\Vert _{L_2(\Pi)}^{2}=\sum_{(i,j)}\pi_{ij}A^2_{ij}$ and Assumption \ref{sampling_1} implies 
$$\Vert A\Vert^{2} _{L_2(\Pi)}\geq p^{-1}\Vert A\Vert^{2} _{2}.$$
\item  Let $\{\epsilon_{ij}\}$ be an i.i.d. Rademacher sequence and $X_{ij}=e_{i}(m_1)e^*_{j}(m_2)$ where $e_{k}(l)$ are the canonical basis vectors in $\mathbb R^l$.
We define
\begin{equation}\label{stoch1}
\Sigma_R=\sum_{(i,j)}\eta_{ij}\epsilon_{ij}X_{ij}\qquad\text{and}\qquad \Sigma=\sum_{(i,j)}\eta_{ij}\xi_{ij}X_{ij}.
\end{equation}

\end{itemize}

\section{The Singular Value Thresholding Algorithm}\label{algorithm}
In this section we introduce an iterative singular value thresholding algorithm and discuss its theoretical properties.  We show that it enjoys strong theoretical guarantees and, unlike one-step thresholding procedures, is well adapted for general non-uniform sampling distributions. 

\subsection{Algorithm}\label{section_algorithm}
 Our algorithm   is based on the \texttt{softImpute} algorithm proposed by Mazumder et al in \cite{Mazumder_Hastie_Tibshirani}. \texttt{SoftImpute} algorithm is inspired by SVD-Impute of Troyanskaya et al \cite{Troyanskaya}.  It alternates between imputing the missing values from a current SVD, and updating the SVD using the data matrix.

\hrulefill

Algorithm 1

\hrulefill

Require : Matrix $Y$, regularization parameter $\lambda$ and $a$, an upper bound on the sup-norm of $M_0$.
\begin{itemize}
\item[1.] $M^{old}=0$
\item[2.]\begin{itemize}
\item[(a)]Repeat\begin{itemize}
\item[(i)] Compute $M^{new}\leftarrow S_{\lambda}\left (Y+(M^{old})_{\bar{\Omega}}\right )$.
\item[(ii)] If $\left \Vert \left (M^{new}-M^{old}\right )_{\bar{\Omega}}\right \Vert<\lambda/3$ and $\left \Vert M^{new}-M^{old}\right \Vert_{\infty}< a$ exit.
\item [(iii)]  Put $M^{old}=\left (M^{old}_{ij}\right )$
 \begin{equation}\label{truncation}
         M^{old}_{ij} = \left\{
    \begin{array}{lll}
                M^{new}_{ij}&  \text{if}\;\;\vert M^{new}_{ij}\vert\leq a \\  \\
  a & \text{if}\; \;M^{new}_{ij}> a \\\\
   -a &\;\text{if}\;\; M^{new}_{ij}< -a.
        \end{array} \right.
        \end{equation}
\end{itemize}
\item[(b)] Assign $\hat M\leftarrow M^{new}$.
\end{itemize}
\item[3.] Output $\hat M$.
\end{itemize}
\hrulefill
\vskip 0.5 cm

This algorithm  repeatedly replaces the missing entries with the current guess,  update the guess by solving
\begin{equation}\label{minimization_problem}
M^{new}\in \underset{M}{\text{minimize}}\; f_{\lambda}(M)=\frac{1}{2}\Vert Y+(M^{old})_{\bar \Omega}-M\Vert^{2}_{2}+\lambda\Vert M\Vert_{*}
\end{equation}
and truncating $M^{new}$.
Let us denote by $(M_{k})_{k\geq 0}$ the sequence of solutions produced by Algorithm 1. We have the following result : 
 
\begin{lemma}\label{lemma_convergency}
For the  successive differences of the  sequence $(M_{k})_{k\geq 0}$ we have that
\begin{equation}\label{convergency_2}
 \left \Vert M^{k+1}-M^{k}\right \Vert_{2}\rightarrow 0\;\text{as}\;k \rightarrow 0
\end{equation}
which implies\label{convergency_infty}
\begin{equation}
\left \Vert \left (M^{k+1}-M^{k}\right )_{\bar{\Omega}}\right \Vert\rightarrow 0 \quad \text{and}\quad \left \Vert M^{k+1}-M^{k}\right \Vert_{\infty}\rightarrow 0\;\text{as}\;k \rightarrow 0.
\end{equation}
\end{lemma}
\subsection{Upper bound on the estimation error}\label{section_upper_bounds}
In this section we derive an upper bound on the estimation error of $\hat M$ produced by Algorithm 1. This bound is non-asymptotic and implies, in particular, that the proposed estimator is minimax optimal. We start by  a general result which  is proven in Appendix \ref{proof-thm-general}. 
\begin{Theorem}\label{thm_general}
Let Assumptions \ref{noise_bounded} and \ref{sampling_1} be satisfied and $\Vert M_0\Vert_{\infty}\leq a$. Assume that  $\lambda\geq 3 \left\Vert \Sigma\right\Vert$.  Then, with probability at least $1-8/d$, 
\begin{equation*} 
\begin{split}
\Vert \hat M-M_0\Vert_{L_{2}(\Pi)}^{2}&\leq C\,p^{-1}\left \{\rank(M_0)\left ( \lambda^{2}+a^{2}\,\left (\bE\left ( \left\Vert \Sigma_R\right\Vert\right )\right )^{2}\right )+{a}^{2}+\log (d)\right \}.
\end{split}
\end{equation*}
 where $d=m_1+m_2$.
\end{Theorem}
Using Assumption \ref{sampling_1}, Theorem \ref{thm_general} implies the following bound on the estimation error measured in normalized Frobenius norm
\begin{Corollary} Under assumptions of Theorem \ref{thm_general} and with probability at least $1-8/d$, 
\begin{equation*} 
\begin{split}
\dfrac{\Vert \hat M-M_0\Vert_{2}^{2}}{m_1m_2}&\leq \dfrac{C}{p^{2}m_1m_2}\left \{\rank(M_0)\left ( \lambda^{2}+a^{2}\,\left (\bE\left ( \left\Vert \Sigma_R\right\Vert\right )\right )^{2}\right )+{a}^{2}+\log (d)\right \}.
\end{split}
\end{equation*}

\end{Corollary}
In order to get a bound in a closed form we need to obtain a suitable upper bounds on $\bE\left ( \left\Vert \Sigma_R\right\Vert\right )$ and, with probability close to $1$, on $\left\Vert \Sigma\right\Vert$. 
 
\begin{lemma}\label{stochastic}
Suppose that $(\xi_{ij})$ are independent and satisfy Assumption \ref{noise_bounded}.   Then, there exists absolute constants $c^{*}, C^{*}>0$  such that, for all $t>0$ with probability at least $1-me^{-t^{2}}$ we have
 \begin{equation}\label{Sigma}
 \left\Vert \Sigma\right \Vert\leq 3\sigma\sqrt{2L}+c^{*}b\,t
 \end{equation} 
 where $L\leq 1$ is defined in \eqref{upper_bound_marginals}.
 
  Moreover, we have
 \begin{equation}\label{expectation_sigma}
 \bE \left\Vert \Sigma_{R}\right\Vert\leq C^{*}\left (\sqrt{L}+\sqrt{\log m}\right ).
 \end{equation}
\end{lemma} 
This Lemma is proven in Appendix \ref{proof-stochastic}. 
 
Taking $t=\sqrt{2\log(d)}$ in Lemma \ref{stochastic}, we get that with probability at least $1-1/d$,
\begin{equation*}
 \left\Vert \Sigma\right \Vert\leq 3\sigma\sqrt{2L}+c^{*}b\,\sqrt{2\log(d)},
\end{equation*}
then, we can choose 
\begin{equation}\label{lambda}
\lambda=3\left (3\sigma\sqrt{2L}+c^{*}b\,\sqrt{2\log(d)}\right ).
\end{equation}
With this choice of $\lambda$ we obtain the following Theorem.
\begin{Theorem}\label{thm3}
Let Assumptions \ref{noise_bounded} and \ref{sampling_1} be satisfied and $\Vert M_0\Vert_{\infty}\leq a$.  Then, with probability at least $1-8/d$, 
\begin{equation*} 
\begin{split}
\Vert \hat M-M_0\Vert_{L_{2}(\Pi)}^{2}&\leq C\,p^{-1}\,\rank(M_0)\left \{\left (a\vee \sigma\right )^{2}L+{a}^{2}\log (m)+b^{2}\log (d)\right \}.
\end{split}
\end{equation*}
and
\begin{equation*} 
\begin{split}
\dfrac{\Vert \hat M-M_0\Vert_{2}^{2}}{m_1m_2}&\leq \dfrac{C\,\rank(M_0)}{p^{2}m_1m_2}\left \{\left (a\vee \sigma\right )^{2}L+{a}^{2}\log (m)+b^{2}\log (d)\right \}.
\end{split}
\end{equation*}
\end{Theorem}

\begin{remark}\label{remark:bound_p}
{\rm Note that $\pi_{ij}\geq p$ yields $L\geq Mp$. Then, the upper bound on the estimation error in the Theorem \ref{thm3} is at least a constant times  $\dfrac{\rank(M_0)}{pm}$. So, in order to get a small estimation error, $p$ should be larger then $\dfrac{\rank(M_0)}{m}$. We denote by $n=\sum_{ij}\pi_{ij}$ the expected number of observations. Condition $p\geq\dfrac{\rank(M_0)}{m} $ implies the following condition on $n$
\begin{equation}\label{cond_n}
n\geq C\,\rank(M_0)\, M.
\end{equation}
When the rank of the matrix $M_0$ is small, this necessary number of observations  is close to the number of degree of freedom of the matrix $M_0$, which is $$(m_1+m_2)\rank(M_0)-\left (\rank(M_0)\right )^{2}.$$
}
\end{remark}
Let us restrict our attention to the non-degenerated case $M_0\not =0$ (we can easily include this case replacing $\rank(M_0)$ by $\rank(M_0)\vee 1$). Assuming that the expected number of observations $n$ is not too small, we can get simpler bound on the estimation error. Suppose that $n>c^{*}m\log(d)$. Then, using $$Lm\geq n\geq c^{*}m\log d$$ we get $L\geq c^{*}\log d$ and we can chose $\lambda$ in the following way \begin{equation}\label{new_lambda}
\lambda=18b\sqrt{2L}.
\end{equation}
With this choice of $\lambda$ we get the following bound on the estimation error
\begin{Corollary}
Let Assumptions \ref{noise_bounded} and \ref{sampling_1} be satisfied and $\Vert M_0\Vert_{\infty}\leq a$. Assume that $n\geq c^{*}m\log(d)$ and $M_0\not =0$. Then, with probability at least $1-8/d$, 
\begin{equation*} 
\begin{split}
\dfrac{\Vert \hat M-M_0\Vert_{2}^{2}}{m_1m_2}&\leq \dfrac{C\,\rank(M_0)\,\left (a\vee b\right )^{2}L}{p^{2}m_1m_2}.
\end{split}
\end{equation*}
 \end{Corollary} 
In order to compare this result with previous results on noisy matrix completion we consider a more restrictive assumption on the sampling distribution. That is, we  assume that this distribution is close to the uniform one:
\begin{Assumption}\label{assumption_uniform}
There exists positives constants $\mu_1$ and $\mu_2$ independent on $m_1$ and $m_2$ and a $0<p<1$ such that for every $(i,j)\in \{1,\dots,m_1\}\times \{1,\dots,m_2\}$ we have

$$\mu_2p\leq \pi_{ij}\leq \mu_1p.$$ 
\end{Assumption}
Under this assumption Theorem \ref{thm_general} yields
\begin{Corollary}\label{corollary_uniforme}
Let Assumptions \ref{noise_bounded} and \ref{assumption_uniform} be satisfied and $\Vert M_0\Vert_{\infty}\leq a$. Assume that $n\geq m\log(d)$ and $\lambda$ given by \eqref{new_lambda}. Then, with probability at least $1-8/d$, 
\begin{equation*} 
\begin{split}
\dfrac{\Vert \hat M-M_0\Vert_{2}^{2}}{m_1m_2}&\leq \dfrac{C\,\rank(M_0)\,\left (a\vee b\right )^{2}}{pm}.
\end{split}
\end{equation*}
\end{Corollary}
\begin{remark}
{\rm Let us compare the bound given by Corollary \ref{corollary_uniforme} with bounds available in the literature. Our model was previously considered by Chatterjee in \cite{Chatterjee_mc} in the case of uniform sampling distribution, that is    $\pi_{ij}=p$ for any $(i,j)\in \{1,\dots,m_1\}\times \{1,\dots,m_2\}$. In \cite{Chatterjee_mc}, Chatterjee introduces a simple estimation procedure, called Universal Singular Value Thresholding which is applied to a number of
questions in low rank matrix estimation,
blockmodels, distance matrix completion, latent space models and etc. For matrix completion problem and under the additional assumption $p\geq n^{-1+\epsilon}$ for some $\epsilon>0$,  the bound obtained in \cite{Chatterjee_mc} is the following one
\begin{equation*} 
\begin{split}
\dfrac{\Vert \hat M-M_0\Vert_{2}^{2}}{m_1m_2}&\leq C\sqrt{\dfrac{\rank(M_0)\,\left (a\vee b\right )^{2}}{pm}}.
\end{split}
\end{equation*}
The rate of convergence given by Corollary \ref{corollary_uniforme} is faster and, as we will see in Section\ref{lower_bounds}, is minimax optimal. Note that the additional assumption $p\geq n^{-1+\epsilon}$ yields the following condition on the expected  number of  observations
\begin{equation}
n>m^{\epsilon}M.
\end{equation}
For low rank matrices, this necessary number of observations is larger than the number of observations required by our method and given by \eqref{cond_n}.

In \cite{Negahban_Wainwright,Koltchinskii_Lounici_Tsybakov, Klopp_general} a closely related set up for matrix completion problem using  the trace regression model was considered. The main difference between these two settings is that in the case of the trace regression the number of observations is not random and each entry may be observed multiple times. In our setting the number of observations is random and each entry is observed at most once. Comparing  with Corollary  \ref{corollary_uniforme} and using $n=pm_1m_2$ we see that bounds obtained in \cite{Negahban_Wainwright,Koltchinskii_Lounici_Tsybakov, Klopp_general} contain an additional logarithmic factor $\log (m_1+m_2)$.}
\end{remark}
\section{Minimax Lower bounds}\label{lower_bounds}
In this section, we prove the minimax lower bound showing that the
rates attained by  our estimator are optimal. The minimax lower bound in a closely related problem  was obtained by Koltchinskii et al in \cite{Koltchinskii_Lounici_Tsybakov}.  We adapt their proof to our set up.

We will denote by $\inf_{\hat{M}}$ the infimum over all the estimators.
 For any $M_0\in \mathbb R^{m_1\times m_2}$, let $\mathbb P_{M_0}$ denote the probability
distribution of the observations $$(\eta_{11}X_{11},\dots,\eta_{m_1m_2}X_{m_1m_2})$$ satisfying \eqref{model}.

For any integer $0\leq r\le
\min(m_1,m_2)$ and any $a>0$, we consider the class of matrices
\begin{equation}\label{Alb_gs}
 \begin{split}
 {\cal A}(r,a)
 &= \left \{ M \in\,\mathbb R^{m_1\times m_2}:\,
 \mathrm{rank}(M)\leq r,\,\Vert M\Vert_{\infty}\leq a,  \right \}.
 \end{split}
 \end{equation}
We will prove the lower bound in the case of the uniform sampling distribution, that is, we suppose that each entry  is observed with the same probability $p$. As it was noted in  Remark \ref{remark:bound_p}, in order to get a small estimation error we need to observe a sufficiently large number of entries, or, equivalently, the  probability  $p$ should be larger then $r/m$. We prove  a lower bound on the estimation risk when this condition is satisfied.
\begin{Theorem}\label{th:lower_bound} Suppose that $m_1,m_2\geq 2$ and $p\geq \frac{r}{m}$.  Fix $a>0$ and integer
$1\leq  r\leq \min(m_1,m_2)$.  Suppose that
the variables $\xi_i$ are i.i.d. Gaussian
${\cal N}(0,\sigma^2)$, $\sigma^2>0$, for $i=1,\dots,n$.  Then, there
exist absolute constants $\beta\in(0,1)$ and $c>0$, such that
\begin{equation*}
\inf_{\hat{M}}
\sup_{\substack{M_0\in\,{\cal A}( r,a)
}}
\mathbb P_{M_0}\left (\dfrac{\Vert \hat M-M_0\Vert_2^{2}}{m_1m_2}> \frac{c\,r\,\left (a\wedge \sigma\right )^{2}}{pm} \right )\ \geq\ \beta.
\end{equation*}
\end{Theorem}
%



\appendix
\section{Proof of Theorem \ref{thm_general}}\label{proof-thm-general}
1. By Lemma 1 in \cite{Mazumder_Hastie_Tibshirani}, $\hat M$ minimizes
\begin{equation*}
f_{\lambda}(M)=\frac{1}{2}\left \Vert Y+(M^{old})_{\bar \Omega}-M\Vert^{2}_{2}+\lambda\Vert M\right \Vert_{*}.
\end{equation*}
Then, using the sub-gradient stationary conditions we have
\begin{equation*}
-\left\langle Y+(M^{old})_{\bar \Omega}-\hat M,\hat M-M_0\right\rangle+\lambda \left\langle \hat V,\hat M-M_0\right\rangle\leq 0
\end{equation*}
where $\hat V\in \partial\Vert \hat M\Vert_{*}$. A simple calculation yields
\begin{equation} \label{thm_gen_1}
\begin{split}
\left\Vert \left (M_0-\hat M\right )_{\Omega}\right \Vert_2^2\leq  \underset{\mathbf {I}}{\underbrace{\left\vert \left\langle \left (Y-M_0\right )_{\Omega},\hat M-M_0\right\rangle\right \vert}}&+\underset{\mathbf {II}}{\underbrace{\left \vert \left\langle \left (M^{old}-\hat M\right )_{\bar \Omega},\hat M-M_0\right\rangle\right \vert}}\\&\hskip 1 cm+ \underset{\mathbf {III}}{\underbrace{\lambda \left\langle \hat V,M_0-\hat M\right\rangle}}.
\end{split}
\end{equation}

2. We estimate each term in \eqref{thm_gen_1} separately. For the first term, we have that $\left (Y-M_0\right )_{\Omega}=\Sigma$ where $\Sigma=\sum_{(i,j)}\eta_{ij}\xi_{ij}X_{ij}$. Then, by the duality between the nuclear and the operator norms, we obtain
\begin{equation}\label{thm_gen_2}
\left\vert \left\langle \left (Y-M_0\right )_{\Omega},\hat M-M_0\right\rangle\right \vert\leq \Vert \Sigma\Vert \Vert \hat M-M_0\Vert_{*}.
\end{equation}
For the second term, using again the duality between the nuclear and the operator norms and the stopping criteria for the Algorithm 1, we obtain
\begin{equation}\label{thm_gen_3}
\begin{split}
\left \vert \left\langle \left (M^{old}-\hat M\right )_{\bar \Omega},\hat M-M_0\right\rangle\right \vert&\leq \left \Vert \left (M^{old}-\hat M\right )_{\bar \Omega}\right  \Vert \left \Vert\hat M-M_0\right \Vert_{*}\\&\hskip 0.5 cm\leq \lambda/3\,\left \Vert \hat M-M_0\right \Vert_{*}.
\end{split}\end{equation}

3. In order to estimate the third term, we use that by monotonicity of subdifferentails of convex functions we have that $\left\langle \hat V-V,\hat M-M_0\right\rangle\geq 0$, for any $ V\in \partial\Vert  M_0\Vert_{*}$. This implies
\begin{equation}\label{thm_general_3}
\left\langle \hat V,M_0-\hat M\right\rangle\leq \left\langle V,M_0-\hat M\right\rangle.
\end{equation} 
Let $P_S$ be the projector on the linear vector subspace $S$ and
 let $S^\bot$ be the orthogonal complement of $S$.
 Let $u_j(A)$ and $v_j(A)$ denote respectively the \textit{left} and \textit{right} orthonormal \textit{singular vectors} of a matrix $A$. $S_1(A)$ is the linear span of $\{u_j(A)\}$, $S_2(A)$ is the linear span of $\{v_j(A)\}$.
 We set
  \begin{equation}\label{projector}
  \begin{split}
  \mathbf P_A^{\bot}(B)=P_{S_1^{\bot}(A)}BP_{S_2^{\bot}(A)}\quad\text{and}\quad\mathbf P_A(B)=B- \mathbf P_A^{\bot}(B).
  \end{split}
  \end{equation}
   Since $\mathbf P_A(B)=P_{S_1^{\bot}(A)}BP_{S_2(A)}+ P_{S_1(A)}B$
  and $\rank(P_{S_i(A)}B)\leq \rank (A)$ we have that
   \begin{equation}\label{thm_gen_9}
  \rank(\mathbf P_A(B))\leq 2\,\rank(A).
  \end{equation}
Note that the subdifferential of the convex function $A\rightarrow  \Vert A \Vert_*$ is the following set of matrices (cf. \cite{watson})
 \begin{equation}\label{subdiff}
 \partial \Vert A \Vert_*=\left\{ \underset{j=1}{\overset{\rank(A)}{\sum}}u_j(A)v_j^{T}(A)+  \mathbf P^{\bot}_{A}(W)\;:\; \Vert  W\Vert\leq 1\right\}.
 \end{equation} 
Inequality \eqref{thm_gen_3} and \eqref{subdiff} imply
\begin{equation}\label{thm_gen_4}
\mathbf {III}\leq \lambda  \left\langle \overset{R}{\sum_{j=1}}u_j(M_0)v_j^{T}(M_0),M_0-\hat M\right\rangle 
+\left\langle  \mathbf P^{\bot}_{M_0}(W),M_0-\hat M\right\rangle.
\end{equation}
Using the fact that $\left \Vert \overset{R}{\sum_{j=1}}u_j(M_0)v_j^{T}(M_0)\right \Vert=1$ and $$\left\langle \overset{R}{\sum_{j=1}}u_j(M_0)v_j^{T}(M_0),M_0-\hat M\right\rangle=\left\langle \overset{R}{\sum_{j=1}}u_j(M_0)v_j^{T}(M_0),\mathbf{P}_{M_0}\left (M_0-\hat M\right )\right\rangle$$
we obtain
\begin{equation}\label{thm_gen_5}
\mathbf {III}\leq \lambda  \left\Vert\mathbf{P}_{M_0}\left (M_0-\hat M\right )\right \Vert_{*}
+\left\langle  \mathbf P^{\bot}_{M_0}(W),M_0-\hat M\right\rangle.
\end{equation}
Now, by the duality between the nuclear and the operator norms, there exists $W$ with $\Vert W\Vert\leq 1$ and such that 
\begin{equation}\label{thm_general_5}
\begin{split}
\left\langle  \mathbf P^{\bot}_{M_0}(W),M_0-\hat M\right\rangle=-\left\langle  W ,\mathbf P^{\bot}_{M_0}\left (\hat M\right )\right\rangle=-\left \Vert \mathbf P^{\bot}_{M_0}\left (\hat M\right )\right \Vert_{*}.
\end{split}\end{equation}
For this particular choice of $W$, \eqref{thm_gen_5} and \eqref{thm_general_5} imply
\begin{equation}\label{thm_gen_6}
\mathbf {III}\leq \lambda  \left (\left\Vert\mathbf{P}_{M_0}\left (M_0-\hat M\right )\right \Vert_{*}
-\left \Vert \mathbf P^{\bot}_{M_0}\left (\hat M\right )\right \Vert_{*}\right ).
\end{equation}
Putting \eqref{thm_gen_2}, \eqref{thm_gen_3}, and \eqref{thm_gen_6} into \eqref{thm_gen_1} and using $\lambda\geq 3\,\Vert \Sigma\Vert$ we obtain
 \begin{equation} 
 \begin{split}
 \left\Vert \left (M_0-\hat M\right )_{\Omega}\right \Vert_2^2\leq  \frac{2\lambda}{3}\,\left \Vert \hat M-M_0\right \Vert_{*}+\lambda  \left (\left\Vert\mathbf{P}_{M_0}\left (M_0-\hat M\right )\right \Vert_{*}
 -\left \Vert \mathbf P^{\bot}_{M_0}\left (\hat M\right )\right \Vert_{*}\right ).
 \end{split}
 \end{equation}

4. The triangle inequality and \eqref{thm_gen_9} lead to

\begin{equation} \label{thm_gen_7}
 \begin{split}
 \left\Vert \left (M_0-\hat M\right )_{\Omega}\right \Vert_2^2\leq  \frac{5\lambda}{3}\,\left\Vert\mathbf{P}_{M_0}\left (M_0-\hat M\right )\right \Vert_{*}\leq \frac{5\lambda\,\sqrt{2\rank(M_0)}}{3}\,\left\Vert M_0-\hat M\right \Vert_{2}
 \end{split}
 \end{equation} 
 and
 \begin{equation} \label{thm_gen_8}
  \begin{split}
  \dfrac{\lambda}{3}\left \Vert \mathbf P^{\bot}_{M_0}\left (\hat M\right )\right \Vert_{*}\leq  \frac{5\lambda}{3}\,\left\Vert\mathbf{P}_{M_0}\left (M_0-\hat M\right )\right \Vert_{*}.
  \end{split}
  \end{equation}
 Inequality \eqref{thm_gen_8} implies 
 \begin{equation*}
  \left \Vert \mathbf P^{\bot}_{M_0}\left (\hat M\right )\right \Vert_{*}\leq  5\,\left\Vert\mathbf{P}_{M_0}\left (M_0-\hat M\right )\right \Vert_{*}
 \end{equation*}
 and
 
   \begin{equation}\label{revision_condition}
   \begin{split}
   \left\Vert \hat M-M_0\right\Vert_{*}\leq 6\left\Vert \mathbf P_{M_0}(\hat M-M_0) \right\Vert_*\leq \sqrt{72\,\rank(M_0)} \,\left\Vert \hat M-M_0\right\Vert_{2}.
   \end{split}
   \end{equation}

5. For a $0<r\leq m$ we consider the following constrain set
\begin{equation}\label{constrain}
\mathcal{C}(r)=\left \{A\in\mathbb{R}^{m_1\times m_2} \,:\, \left\Vert A\right\Vert_{\infty}=1, \left\Vert A\right\Vert_{L_2(\Pi)}^{2}\geq \dfrac{\log(d)}{0.0006\,\log\left (6/5\right )\,p}, \left\Vert A\right\Vert_{*}\leq \sqrt{r} \left\Vert A\right\Vert_{2}\right \}.
\end{equation}
Note that the condition $\left\Vert A\right\Vert_{*}\leq \sqrt{r} \left\Vert A\right\Vert_{2}$ is satisfied if $\rank(A)\leq r$.

We have the following result for matrices in  $\mathcal{C}(r)$. Its proof is given in Appendix \ref{proof-thm1}.
\begin{lemma}\label{thm1}
For all $A\in \mathcal{C}(r)$
$$\left\Vert A_{\Omega}\right\Vert^{2}_{2}\geq \frac{\Vert A\Vert _{L_2(\Pi)}^{2}}{2}-44\,p^{-1}\left [r\left (\bE\left ( \left\Vert \Sigma_R\right\Vert\right )\right )^{2}+18\right ]$$
with probability at least $1-8/d$.
\end{lemma}

Note that condition $\left \Vert \hat M-M^{old}\right \Vert_{\infty}< a$ and $\left \Vert M^{old}\right \Vert_{\infty}< a$ imply $$\left\Vert \hat M-M_0\right\Vert_{\infty}\leq 3a.$$
 We now consider two cases, depending on whether the matrix $\dfrac{\left (\hat M-M_0\right )}{3a}$ belongs to the set $ \mathcal{C}\left (72\,\rank(M_0)\right )$ or not.

\textbf{Case 1}: Suppose first that $\left\Vert \hat M-M_0\right\Vert_{L_2(\Pi)}^{2} < \dfrac{\log(d)}{0.0006\,\log\left (6/5\right )\,p}$, then the statement  of the Theorem \ref{thm_general} is true.

\textbf{Case 2}: It remains to consider the case  $\left\Vert \hat M-M_0\right\Vert_{L_2(\Pi)}^{2} \geq \dfrac{\log(d)}{0.0006\,\log\left (6/5\right )\,p}$. Then \eqref{revision_condition} implies that $\dfrac{1}{3a}\left (\hat M-M_0\right )\in \mathcal{C}\left (72\,\rank(M_0)\right )$ and we can apply Lemma \ref{thm1}. From Lemma \ref{thm1} and \eqref{thm_gen_7} we obtain that 
with probability at least $1-8/d$ one has
\begin{equation*} 
 \begin{split}
 \dfrac{1}{2}\Vert \hat M-M_0\Vert _{L_2(\Pi)}^{2}&\leq \frac{5\lambda\,\sqrt{2\rank(M_0)}}{3}\,\left\Vert M_0-\hat M\right \Vert_{2}\\&+369\,a^{2}\,p^{-1}\left [72\,\rank(M_0)\left (\bE\left ( \left\Vert \Sigma_R\right\Vert\right )\right )^{2} +18\right ]
 \\&\leq 6\lambda^{2}\,p^{-1}\rank(M_0)+\frac{p}{4}\left\Vert \hat M-M_0 \right\Vert^{2}_2 \\&\hskip 0.25 cm+369\,a^{2}\,p^{-1}\left [72\,\rank(M_0)\left (\bE\left ( \left\Vert \Sigma_R\right\Vert\right )\right )^{2} +18\right ].
 \end{split}
 \end{equation*}
Now \eqref{ass1} imply that, there exist numerical constants $C$ such that
\begin{equation*} 
\begin{split}
\Vert \hat M-M_0\Vert_{L_2(\Pi)}^{2}&\leq C\,p^{-1}\left \{\rank(M_0)\left ( \lambda^{2}\,+a^{2}\,\left (\bE\left ( \left\Vert \Sigma_R\right\Vert\right )\right )^{2}\right )+a^{2}\right \},
\end{split}
\end{equation*}
which leads to the statement of the Theorem \ref{thm_general}.
\section{Proof of Theorem \ref{th:lower_bound}}\label{proof_th:lower_bound}
We adopt the proof of Theorem 5 in \cite{Koltchinskii_Lounici_Tsybakov} to our setting. Assume w.l.o.g. that $m_1\geq m_2$. For a $\gamma\leq 1$, define
$$
\tilde{ \mathcal{L}} \, =\left \{ \tilde{L} = (l_{ij})\in\R^{m_1\times r}:
l_{ij}\in\left \{0, \gamma(\sigma \wedge a) \left (\frac{
 r}{p\,m}\right )^{1/2}\right \}\,, \forall 1\leq i \leq m_1,\, 1\leq j\leq
r\right \},
$$
and consider the associated set of block matrices
$$
\mathcal{A} \ =\ \Big\{
L=(\begin{array}{c|c|c|c}\tilde{L}&\cdots&\tilde{L}&O
\end{array})\in\R^{m_1\times m_2}: \tilde{L}\in \tilde{\mathcal{L}}\Big\},
$$
where $O$ denotes the $m_1\times (m_2-r\lfloor m_2/(2r) \rfloor )$ zero
matrix, and $\lfloor x \rfloor$ is the integer part of $x$.

\begin{remark}
In the case $m_1< m_2$, we only need to change the construction of the low rank component of the test set. We first build a matrix $\tilde L = \left(\begin{array}{c|c}\bar L&O\\\end{array}\right) \in \mathbb R^{r \times m_2} $ where $\bar L \in \mathbb \mathbb R^{r \times (m_2/2)}$ with entries in $\left\{0, \gamma (\sigma \wedge a) \left(\frac{ r}{p\,m} \right)^{1/2}\right\}$ and, then, we replicate this matrix to obtain a block matrix $L$ of size $m_1 \times m_2$ 
$$
L=\left(
\begin{array}{c}
\tilde{L}\\
\hline\\
\vdots\\
\hline\\
\tilde{L}\\
\hline\\
O
\end{array}
\right).
$$
\end{remark}

By construction, any element of $\mathcal{A}$ as well
as the difference of any two elements of $\mathcal{A}$ has  rank at most $r$. In addition, condition $p\geq \frac{r}{m}$ implies that the entries of any matrix in
$\mathcal{A}$ take values in $[0,a]$. Thus,
$\mathcal{A}\subset {\cal A}(r,a)$. 

   The Varshamov-Gilbert bound (cf. Lemma 2.9 in \cite{tsy_09}) guarantees the existence of a subset $\mathcal A^0\subset{\mathcal{A}}$ with
cardinality $\mathrm{Card}(\mathcal A^0) \geq 2^{(rM)/8}+1$ containing the
zero $m_1\times m_2$ matrix ${\bf 0}$ and such that, for any two
distinct elements $A_1$ and $A_2$ of $\mathcal A^0$,
\begin{equation}\label{lower_2}
\Arrowvert A_1-A_2\Arrowvert_{2}^2  \geq \frac{Mr}{8}
\left(\gamma^2(\sigma \wedge a)^2 \frac{r }{p\,m} \right)
\left\lfloor \frac{m}{r}\right\rfloor \geq
\frac{\gamma^2}{16}(\sigma \wedge a)^2 \,m_1m_2\frac{r}{p\,m}\,.
\end{equation}

Using that, conditionally on $X_i$, the distributions of $\xi_i$ are
Gaussian, we get that, for any $A\in \mathcal A_0$, the Kullback-Leibler
divergence $K\big(\mathbb P_{{\bf 0}},\mathbb P_{A}\big)$ between $\mathbb  P_{{\bf 0}}$
and $\mathbb  P_{A}$ satisfies
\begin{equation}\label{KLdiv}
K\big(\mathbb P_{{\bf 0}},\mathbb P_{A}\big)\ =\
\frac{1}{2\sigma^2}\|A\|_{L_2(\Pi)}^2 \leq
\frac{\gamma^2\,Mr}{2}.
\end{equation}
From (\ref{KLdiv}) we deduce that the condition
\begin{equation}\label{eq: condition C}
\frac{1}{\mathrm{Card}({\cal A}^0)-1} \sum_{A\in{\cal A}^0}K(\mathbb  P_{\bf
0},\mathbb P_{A})\ \leq\ \alpha \log \big(\mathrm{Card}({\cal A}^0)-1\big)
\end{equation}
is satisfied for any $\alpha>0$ if $ \gamma>0$ is chosen as a
sufficiently small numerical constant depending on $\alpha$.  In view
of (\ref{lower_2}) and (\ref{eq: condition C}) and using  the application of Theorem 2.5 in \cite{tsy_09} implies
\begin{equation}\label{lower_1}
\inf_{\hat{M}}
\sup_{\substack{M_0\in\,{\cal A}( r,a)
}}
\mathbb P\left (\dfrac{\Vert \hat M-M_0\Vert_2^{2}}{m_1m_2}> \dfrac{C(\sigma \wedge a)^2\,r}{pm} \right )\ \geq\ \beta
\end{equation}
for some  absolute constants $\beta\in(0,1)$, which implies the statement of Theorem \ref{th:lower_bound}.

\section{Proof of Lemma \ref{thm1}}\label{proof-thm1}
This proof is close to the proof of Lemma 12 in \cite{Klopp_general}. Set $$\mathcal{E}=44\,p^{-1}\left [r\left (\bE\left ( \left\Vert \Sigma_R\right\Vert\right )\right )^{2} +18\right ].$$
We will show that the probability of the following ``bad'' event is small
\begin{equation*}
\mathcal{B}=\left \{\exists\,A\in \mathcal{C}(r)\,\text{such that}\,\left \vert  \left\Vert A_{\Omega} \right\Vert_{2}^{2}-\Vert A\Vert _{L_2(\Pi)}^{2}\right \vert> \dfrac{1}{2}\Vert A\Vert _{L_2(\Pi)}^{2}+ \mathcal{E}\right \}.
\end{equation*}
Note that $\mathcal{B}$ contains the complement of the event that we are interested in.

In order to estimate the probability of $\mathcal{B}$ we use a standard peeling argument. Let $\nu=\dfrac{\log(d)}{0.0006\,\log\left (6/5\right )\,p}$ and $\alpha=\dfrac{6}{5}$. For $l\in\mathbb N$ set $$S_l=\left \{A\in \mathcal{C}(r)\,:\,\alpha^{l-1}\nu \leq \Vert A\Vert _{L_2(\Pi)}^{2}\leq \alpha^{l}\nu\right \}.$$ If the event $\mathcal{B}$ holds for some matrix $A\in \mathcal{C}(r)$, then $A$ belongs to some $S_l$ and 
\begin{equation}\label{Bl}
\begin{split}
\left \vert \left\Vert A_{\Omega} \right\Vert_{2}^{2}-\Vert A\Vert _{L_2(\Pi)}^{2}\right \vert&> \dfrac{1}{2}\Vert A\Vert _{L_2(\Pi)}^{2}+ \mathcal{E}\\&> \dfrac{1}{2}\alpha^{l-1}\nu+ \mathcal{E}\\&
= \dfrac{5}{12}\alpha^{l}\nu+ \mathcal{E}.
\end{split}
\end{equation}
For $T>\nu$ consider the following set of matrices
$$\mathcal{C}(r,T)=\left \{A\in\mathcal{C}(r) \,:\,  \left\Vert A\right\Vert_{L_2(\Pi)}^{2}\leq T \right \}
$$ and the following event 
$$\mathcal{B}_l=\left \{\exists\,A\in \mathcal{C}(r,\alpha^{l}\nu)\,:\,\left \vert \left\Vert A_\Omega \right\Vert_{2}^{2}-\Vert A\Vert _{L_2(\Pi)}^{2}\right \vert> \dfrac{5}{12}\alpha^{l}\nu+ \mathcal{E}\right \}.$$
Note that $A\in S_l$ implies that $A\in \mathcal{C}(r,\alpha^{l}\nu)$. Then \eqref{Bl} implies that $\mathcal{B}_l$ holds and we get $\mathcal{B}\subset\cup \,\mathcal{B}_l$. Thus, it is enough to estimate the probability of the simpler event $\mathcal{B}_l$ and then apply the union bound. Such an estimation is given by the following lemma. Its proof is given in Appendix \ref{pl1}. Let
$$Z_T=\underset{A\in \mathcal{C}(r,T)}{\sup}\left \vert \left\Vert A_\Omega \right\Vert_{2}^{2}-\Vert A\Vert _{L_2(\Pi)}^{2}\right \vert.$$ 
\begin{lemma}\label{l1}
We have that  \begin{equation*}
 \mathbb{P}\left ( Z_{T} \geq \frac{5}{12}T+44\,p^{-1}\left [r\left (\bE\left ( \left\Vert \Sigma_R\right\Vert\right )\right )^{2}+18\right ]\right )\leq 4e^{-c_1\,p\,T}
 \end{equation*}
 with $c_1\geq 0.0006$.
 \end{lemma}
Lemma \ref{l1} implies that $\mathbb P\left (\mathcal{B}_l\right )\leq 4\exp(-c_1\,p\,\alpha^{l}\nu)$. Using the union bound we obtain
\begin{equation*}
\begin{split}
\mathbb P\left (\mathcal{B}\right )&\leq \underset{l=1}{\overset{\infty}{\Sigma}}\mathbb P\left (\mathcal{B}_l\right )\\&\leq 4\underset{l=1}{\overset{\infty}{\Sigma}}\exp(-c_1\,p\,\alpha^{l}\nu)\\&\leq 4\underset{l=1}{\overset{\infty}{\Sigma}}\exp\left (-c_1\,p\,\nu\,\log(\alpha)\,l\right )
\end{split}
\end{equation*} 
where we used $e^{x}\geq x$. We finally compute for $\nu=\dfrac{\log(d)}{0.0006\,p\,\log\left (6/5\right )}$
\begin{equation*}
\mathbb P\left (\mathcal{B}\right )\leq \dfrac{4\exp\left (-c_1\,p\,\nu\,\log(\alpha)\right )}{1-\exp\left (-c_1\,p\,\nu\,\log(\alpha)\right )}=\dfrac{4\exp\left (-\log(d)\right )}{1-\exp\left (-\log(d)\right )}.
\end{equation*} 
This completes the proof of Lemma \ref{thm1}.\\
\section{Proof of Lemma \ref{l1}}\label{pl1}
 
We will start by showing that $Z_T$ concentrates around its expectation and then we will upper bound the expectation. Recall that by definition, $$Z_T=\underset{A\in \mathcal{C}(r,T)}{\sup}\left \vert \sum _{(i,j)} \eta_{ij}A^{2}_{ij}-\bE\left (\sum _{(i,j)} \eta_{ij}A^{2}_{ij}\right )\right \vert.$$
We use the following Talagrand's concentration inequality :
\begin{Theorem}\label{talagrand}
Suppose that $f\,:\,[-1,1]^{N}\rightarrow \mathbb{R}$ is a convex Lipschitz function with Lipschitz constant $L$. Let $\Xi_1,\dots \Xi_N$ be independent random variables taking value in $[-1,1]$. Let $Z:\,=f(\Xi_1,\dots, \Xi_n)$. Then for any $t\geq 0$,
$$\mathbb{P}\left (\left \vert Z-\mathbb{E}(Z)\right \vert\geq 16L+t\right )\leq 4e^{-t^{2}/2L^{2}}.$$
\end{Theorem}
For a proof see \cite{talagrand1996} and \cite{Chatterjee_mc}. Let $f(x_{11},\dots,x_{m_1m_2}):\,= \underset{A\in \mathcal{C}(r,T)}{\sup}\left \vert \sum _{(i,j)} \left (x_{ij}-p_{ij}\right )A^{2}_{ij}\right \vert.$   It is easy to see that $f(x_{11},\dots,x_{m_1m_2})$ is a Lipschitz function with Lipschitz constant $L=\sqrt{p^{-1}T}$. Indeed,
 \begin{equation*}
 \begin{split}
 \left \vert f(x_{11},\dots,x_{m_1m_2})-f(z_{11},\dots,z_{m_1m_2})\right \vert& \\&\hskip -3 cm =\left \vert \underset{A\in \mathcal{C}(r,T)}{\sup}\left \vert \sum _{(i,j)} \left (x_{ij}-p_{ij}\right )A^{2}_{ij}\right \vert-\underset{A\in \mathcal{C}(r,T)}{\sup}\left \vert \sum _{(i,j)} \left (z_{ij}-p_{ij}\right )A^{2}_{ij}\right \vert\right \vert\\
  &\hskip -2.5 cm\leq  \underset{A\in \mathcal{C}(r,T)}{\sup}\left \vert  \left \vert\sum _{(i,j)} \left (x_{ij}-p_{ij}\right )A^{2}_{ij}\right \vert- \left \vert\sum _{(i,j)} \left (z_{ij}-p_{ij}\right )A^{2}_{ij}\right \vert\right \vert\\
 &\hskip -2 cm \leq  \underset{A\in \mathcal{C}(r,T)}{\sup}\left \vert  \sum _{(i,j)} \left (x_{ij}-p_{ij}\right )A^{2}_{ij}- \sum _{(i,j)} \left (z_{ij}-p_{ij}\right )A^{2}_{ij}\right \vert\\& \hskip -1.5 cm\leq \underset{A\in \mathcal{C}(r,T)}{\sup}\left \vert  \sum _{(i,j)} \left (x_{ij}-z_{ij}\right )A^{2}_{ij}\right \vert
 \\&\hskip -1 cm \leq \underset{A\in \mathcal{C}(r,T)}{\sup}\sqrt{\sum _{(i,j)} \pi^{-1}_{ij}\left (x_{ij}-z_{ij}\right )^{2}} \sqrt{\sum _{(i,j)} \pi_{ij}A^{4}_{ij}}\\&\hskip -0.5 cm \leq\sqrt{p^{-1}}\underset{A\in \mathcal{C}(r,T)}{\sup}\sqrt{\sum _{(i,j)} \left (x_{ij}-z_{ij}\right )^{2}}\sqrt{\sum _{(i,j)} \pi_{ij}A^{2}_{ij}}
 \\&\hskip 0.5 cm \leq\sqrt{p^{-1}T}\sqrt{\sum _{(i,j)} \left (x_{ij}-z_{ij}\right )^{2}}
 \end{split}
 \end{equation*}
 where we used $\left \vert \vert a\vert- \vert b \vert \right \vert\leq \left \vert a-b\right \vert$, $\Vert A\Vert _{\infty}\leq 1$ and $\left\Vert A\right\Vert_{L_2(\Pi)}^{2}\leq T$. Now, Theorem \ref{talagrand} and $2\sqrt{p^{-1}T}\leq T+p^{-1} $ imply
 $$\mathbb{P}\left ( Z_{T}\geq \mathbb{E}(Z_{T})+768\,p^{-1}+\frac{1}{12}T+t\right )\leq 4e^{-t^{2}p/2\,T}.$$
 Taking $t=\frac{1}{9}\left (\frac{1}{3}T\right )$ we get
 \begin{equation}\label{concentration_bound}
 \mathbb{P}\left ( Z_{T} \geq \mathbb{E}(Z_{T})+ 768\,p^{-1}+\frac{1}{9}\left (\frac{5}{12}T\right )\right )\leq 4e^{-c_1\,p\,T}
 \end{equation}
 with $c_1\geq 0.0006$.
  
Next we bound the expectation $\bE\left ( Z_T\right )$. Using a standard symmetrization argument (see e.g. \cite{Ledoux_book}) we obtain 
\begin{equation*}
\begin{split}
\bE \left ( Z_T\right )&= \bE\left (\underset{A\in \mathcal{C}(r,T)}{\sup}\left \vert \sum_{(i,j)} \eta_{ij}A_{ij}^{2}-\bE\left (\eta_{ij}A_{ij}^{2}\right )\right \vert\right )\\&\leq 2\bE\left (\underset{A\in \mathcal{C}(r,T)}{\sup}\left \vert\sum_{(i,j)} \epsilon_{ij}\eta_{ij}A_{ij}^{2} \right \vert\right )
\end{split}\end{equation*}
where $\{\epsilon_{ij}\}$ is an i.i.d. Rademacher sequence. Then, the contraction inequality (see e.g. \cite[Theorem 2.2]{Koltchinskii_st_flour}) yields
\begin{equation*}
\bE \left ( Z_T\right )\leq 8\bE\left (\underset{A\in \mathcal{C}(r,T)}{\sup}\left \vert\sum_{(i,j)}\epsilon_{ij}\eta_{ij}A_{ij}\right \vert\right )=8\,\bE\left (\underset{A\in \mathcal{C}(r,T)}{\sup}\left \vert \left\langle \Sigma_R,A\right\rangle\right \vert\right )
\end{equation*}
where $\Sigma_R=\sum_{(i,j)}\epsilon_{ij}\eta_{ij}X_{ij}$. For $A\in \mathcal{C}(r,T)$ we have that 
\begin{equation*}
\begin{split}
\left\Vert A \right\Vert_* &\leq \sqrt{r}\left\Vert A\right\Vert_2\\&\leq  \sqrt{r\,p^{-1}}\left\Vert A\right\Vert_{L_2(\Pi)}\\&\leq \sqrt{r\,p^{-1}\,T}
\end{split}
\end{equation*}
where we have used \eqref{ass1}. Then, by the duality between nuclear and operator norms, we compute 
\begin{equation*}
\bE \left ( Z_T\right )\leq 8\bE\left (\underset{\left\Vert A \right\Vert_*\leq \sqrt{r\,p^{-1}\,T}}{\sup}\left \vert \left\langle \Sigma_R,A\right\rangle\right \vert\right )\leq 8\,\sqrt{r\,p^{-1}\,T}\,\bE\left ( \left\Vert \Sigma_R\right\Vert\right ).
\end{equation*}
Finally, using $$\dfrac{1}{9}\left (\dfrac{5}{12}T\right )+8\sqrt{r\,p^{-1}\,T}\,\bE\left ( \left\Vert \Sigma_R\right\Vert\right )\leq \left (\dfrac{1}{9}+\dfrac{8}{9}\right ) \dfrac{5}{12}T+44\,r\,p^{-1}\left (\bE\left ( \left\Vert \Sigma_R\right\Vert\right )\right )^{2}$$
and the concentration bound \eqref{concentration_bound} we obtain that 
 \begin{equation*}
 \mathbb{P}\left ( Z_{T} \geq \frac{5}{12}T+ 44\,p^{-1}\left [r\left (\bE\left ( \left\Vert \Sigma_R\right\Vert\right )\right )^{2}+18\right ]\right )\leq 4e^{-c_1\,p\,T}
 \end{equation*}
 with $c_1\geq 0.0006$ as stated.
\section{Proof of Lemma \ref{lemma_convergency}}\label{pl2}
It is easy to see that
\begin{equation*}
\left \Vert \left (M_{k+1}-M_{k}\right )_{\bar{\Omega}}\right \Vert\leq \left \Vert \left (M_{k+1}-M_{k}\right )_{\bar{\Omega}}\right \Vert_{2}\leq \left \Vert M_{k+1}-M_{k}\right \Vert_{2}
\end{equation*}
and
\begin{equation*}
\left \Vert M_{k+1}-M_{k}\right \Vert_{\infty}\leq \left \Vert M_{k+1}-M_{k}\right \Vert_{2}.
\end{equation*}
Thus, it is enough to show \eqref{convergency_2}. The proof of \eqref{convergency_2} is close  to the proof of  Lemma 4 in \cite{Mazumder_Hastie_Tibshirani}.

Let us denote for by $\tilde M_{k}$  the  solutions produced by Algorithm 1 after soft-thresholding step and before truncating step \eqref{truncation}. We have that
\begin{equation}\label{lemma3_1}
\Vert M^{k+1}-M^{k}\Vert_{2}\leq \Vert \tilde M^{k+1}-\tilde M^{k}\Vert_{2}\leq \Vert (M^{k}-M^{k-1})_{\bar \Omega}\Vert_{2}\leq \Vert M^{k}-M^{k-1}\Vert_{2}
\end{equation}
where in the second inequality we used the following result (see, for example, Lemma 3 in \cite{Mazumder_Hastie_Tibshirani})
\begin{proposition}
The soft-thresholding operator $S_{\lambda}(\cdot)$ satisfies the following: for any $W_1,W_2$
$$\Vert S_{\lambda}(W_1)-S_{\lambda}(W_2)\Vert _{2}\leq \Vert W_1-W_2\Vert_{2}.$$
\end{proposition}

 The inequality \eqref{lemma3_1} implies that the sequence $\{\Vert M^{k}-M^{k-1}\Vert_{2}\}_{k\geq 1}$  converges. It remains to show that it converges to zero. Note that the inequalities \eqref{lemma3_1} imply that
$$\Vert M^{k}-M^{k-1}\Vert^{2}_{2}-\Vert (M^{k+1}-M^{k})_{\bar \Omega}\Vert^{2}_{2}=\Vert (M^{k+1}-M^{k})_{\Omega}\Vert^{2}_{2}\rightarrow 0.$$
So, we only need to show that $\Vert (M^{k+1}-M^{k})_{\bar \Omega}\Vert_{2}\rightarrow 0$.

 We put $$Q(A,B)=\frac{1}{2}\Vert (Y-B)_{\Omega}\Vert^{2}_{2}+\frac{1}{2}\Vert (A-B)_{\bar \Omega}\Vert^{2}_{2}+\lambda\Vert B\Vert_{*}.$$ Note that \eqref{minimization_problem} implies 
\begin{equation}\label{lemma3_2}
\begin{split}
Q(M^{k},\tilde M^{k})&\geq Q(M^{k},\tilde M^{k+1})\\&=\frac{1}{2}\Vert (Y-\tilde M^{k+1})_{\Omega}\Vert^{2}_{2}+\frac{1}{2}\Vert (M^{k}-\tilde M^{k+1})_{\bar \Omega}\Vert^{2}_{2}+\lambda\Vert \tilde M^{k+1}\Vert_{*}
\\&\geq \frac{1}{2}\Vert (Y-\tilde M^{k+1})_{\Omega}\Vert^{2}_{2}+\frac{1}{2}\Vert (M^{k+1}-\tilde M^{k+1})_{\bar \Omega}\Vert^{2}_{2}+\lambda\Vert \tilde M^{k+1}\Vert_{*}\\&=Q( M^{k+1},\tilde M^{k+1})
\end{split}
\end{equation}
 where in the last inequality we used that
 \begin{equation}\label{lemma3_5}
          M^{k+1}_{ij} = \left\{
     \begin{array}{lll}
                 \tilde M^{k+1}_{ij}&  \text{if}\;\;\vert \tilde M^{k+1}_{ij}\vert\leq a \\  \\
   a & \text{if}\; \;\tilde M^{k+1}_{ij}> a \\\\
    -a &\;\text{if}\;\; \tilde M^{k+1}_{ij}< -a.
         \end{array} \right.
         \end{equation}
 The inequality \eqref{lemma3_2} shows that the sequence $\{Q(M^{k},\tilde M^{k})\}_{k\geq 1}$ converges. This and \eqref{lemma3_2} yield 
 \begin{equation}\label{lemma3_3}
 \begin{split}
 Q(M^{k},\tilde M^{k+1})-Q( M^{k+1},\tilde M^{k+1})\\&\hskip -2 cm=\frac{1}{2}\Vert (M^{k}-\tilde M^{k+1})_{\bar \Omega}\Vert^{2}_{2}-\frac{1}{2}\Vert (M^{k+1}-\tilde M^{k+1})_{\bar \Omega}\Vert^{2}_{2}\rightarrow 0. \end{split}
 \end{equation}
 Now, it is easy to see that 
 \begin{equation}\label{lemma3_4}
 \Vert (M^{k}-\tilde M^{k+1})_{\bar \Omega}\Vert^{2}_{2}-\Vert (M^{k+1}-\tilde M^{k+1})_{\bar \Omega}\Vert^{2}_{2}\geq \Vert (M^{k}- M^{k+1})_{\bar \Omega}\Vert^{2}_{2}.
 \end{equation}
 Indeed, for $(i,j)$ in $\bar \Omega$ such that $M_{ij}^{k+1}=\tilde M_{ij}^{k+1}$ we have that  $$\left (M_{ij}^{k}-\tilde M_{ij}^{k+1}\right) ^{2}-\left (M_{ij}^{k+1}-\tilde M_{ij}^{k+1}\right) ^{2}= \left (M_{ij}^{k}-M_{ij}^{k+1}\right) ^{2}$$
 and for $(i,j)$ in $\bar \Omega$ such that $ M_{ij}^{k+1}\not = M_{ij}^{k+1}$ 
 we have that  $$\left (M_{ij}^{k}-\tilde M_{ij}^{k+1}\right) ^{2}-\left (M_{ij}^{k+1}-\tilde M_{ij}^{k+1}\right) ^{2}\geq \left (M_{ij}^{k}-M_{ij}^{k+1}\right) ^{2}$$
 where we used \eqref{lemma3_5}.
 Now \eqref{lemma3_3} together with \eqref{lemma3_4} imply \eqref{convergency_2}
 which completes the proof  of Lemma \ref{lemma_convergency}.

 \section{Proof of Lemma \ref{stochastic}}\label{proof-stochastic}
 In order to prove \eqref{Sigma}, we use the following remarkable bound on the spectral norms of random matrices. It is obtained by extension to rectangular matrices via self-adjoint dilation of Corollary 3.12 and Remark 3.13 in \cite{Bandeira}  (cf., Section 3.1 in \cite{Bandeira}).
  \begin{proposition}[\cite{Bandeira}]\label{pr1}
  Let $A$ be the $m_1\times m_2$ rectangular matrix whose entries $A_{ij}$ are independent centered bounded random variables. Then, for any $0<\epsilon\leq 1/2$ there exists a universal constant $c_{\epsilon}$ such that, for every $t\geq 0$
  $$ \mathbb P\left \{ \left \Vert A\right\Vert\geq (1+\epsilon)2\sqrt{2}(\sigma_1\vee \sigma_2)+t\right \}\leq (m_1\wedge m_2)\exp\left (\frac{-t^{2}}{c_{\epsilon}\sigma^{2}_*}\right )$$ 
  where we have defined
  $$\sigma_{1}=\underset{i}{\max}\sqrt{\sum_{j}\mathbb E[A_{ij}^{2}]}$$
   $$\sigma_{2}=\underset{j}{\max}\sqrt{\sum_{i}\mathbb E[A_{ij}^{2}]}$$
    $$\sigma_{*}=\underset{ij}{\max}\vert A_{ij}\vert.$$
     \end{proposition}
     
     We apply Proposition \ref{pr1} to $\Sigma=\sum_{(i,j)}\eta_{ij}\xi_{ij}X_{ij}$. We compute
     \begin{equation*}
     \begin{split}
     \sigma_1=\underset{i}{\max}\sqrt{\sum_{j}\mathbb E[\eta_{ij}^{2}\xi_{ij}^{2}]}=
     \sigma \underset{i}{\max}\sqrt{\pi_{i\cdot}}\quad \text{and}\quad \sigma_2=
           \sigma \underset{j}{\max}\sqrt{\pi_{\cdot j}}.
     \end{split}
     \end{equation*}
 Bound \eqref{upper_bound_marginals} implies that $\sigma_1\vee \sigma_2\leq \sigma \sqrt{L}$. On the other hand, Assumption \ref{noise_bounded} implies $\underset{ij}{\max}\vert \eta_{ij}\xi_{ij}\vert\leq b$. Now, taking in Proposition \ref{pr1} $\epsilon=1/2$ we get \eqref{Sigma}.
 
 In order to prove \eqref{expectation_sigma} we use the following result
  \begin{proposition}[Corollary 3.3 in \cite{Bandeira}]\label{pr2}
   Let $A$ be the $m_1\times m_2$ rectangular matrix with $A_{ij}$ independent centered bounded random  variables. Then,  there exists a universal constant $C^{*}$ such that, 
   $$ \mathbb E \left \Vert A\right\Vert\leq C^{*}\left \{\sigma_1\vee \sigma_2+\sigma_{*}\sqrt{\log(m_1\wedge m_2)}\right \}$$ 
   where $\sigma_1, \sigma_2, \sigma_*$ are defined in Proposition \ref{pr1}.
      \end{proposition}
      We apply Proposition \ref{pr2} to $\Sigma_R=\sum_{(i,j)}\eta_{ij}\epsilon_{ij}X_{ij}$ where $\{\epsilon_{ij}\}$ is i.i.d. Rademacher sequence. We have that $\sigma_1\vee \sigma_2\leq \sqrt{L}$ and $\sigma_{*}\leq 1$, then Proposition \ref{pr2} implies \eqref{expectation_sigma}.




\end{document}